\newif\ifpdf
\ifx\pdfoutput\undefined
\pdffalse 
\else
\pdfoutput=1 
\pdftrue \fi

\newif\iffinal
\finalfalse	
\finaltrue 

\documentclass[reqno,twoside,11pt]{amsart}
\iffinal\else\usepackage[notref,notcite]{showkeys}\fi
\usepackage{amsmath}
\usepackage{amsfonts}
\usepackage{amssymb}
\usepackage{verbatim}
\usepackage{epsfig}
\usepackage{color}

\IfFileExists{epsf.def}{\input epsf.def}{\usepackage{epsf}}

\usepackage{mathpazo}

\DeclareFontFamily{OT1}{eusb}{} \DeclareFontShape{OT1}{eusb}{m}{n} {<5> <6> <7> <8> <9> <10> <11> <12> <14.4> eusb10}{}
\DeclareMathAlphabet{\eusb}{OT1}{eusb}{m}{n}

\DeclareFontFamily{OT1}{eusm}{} \DeclareFontShape{OT1}{eusm}{m}{n} {<5> <6> <7> <8> <9> <10> <11> <12> <14.4> eusm10}{}
\DeclareMathAlphabet{\eusm}{OT1}{eusm}{m}{n}

\DeclareFontFamily{OT1}{eufm}{} \DeclareFontShape{OT1}{eufm}{m}{n} {<5> <6> <7> <8> <9> <10> <11> <12> <14.4> eufm10}{}
\DeclareMathAlphabet{\mathfrak}{OT1}{eufm}{m}{n}

\DeclareFontFamily{OT1}{fraktura}{}
\DeclareFontShape{OT1}{fraktura}{m}{n} {<5> <6> <7> <8> <9> <10> <11> <12> <13> <14.4> [1.1] eufm10}{}
\DeclareMathAlphabet{\fraktura}{OT1}{fraktura}{m}{n}

\DeclareFontFamily{OT1}{cmfi}{} \DeclareFontShape{OT1}{cmfi}{m}{n} {<5> <6> <7> <8> <9> <10> <11> <12> <13> <14.4> [0.9] cmfi10}{}
\DeclareMathAlphabet{\cmfi}{OT1}{cmfi}{b}{n}

\DeclareFontFamily{OT1}{cmss}{} \DeclareFontShape{OT1}{cmss}{m}{n} {<5> <6> <7> <8> <9> <10> <11> <12> <13> <14.4> cmss10}{}
\DeclareMathAlphabet{\cmss}{OT1}{cmss}{m}{n}

\newtheoremstyle{thm}{1.5ex}{1.5ex}{\itshape\rmfamily}{} {\bfseries\rmfamily}{}{2ex}{}

\newtheoremstyle{def}{1.5ex}{1.5ex}{\slshape\rmfamily}{} {\bfseries\rmfamily}{}{2ex}{}

\newtheoremstyle{rem}{1.3ex}{1.3ex}{\rmfamily}{} {\itshape}
{} {1.5ex}{}

\newenvironment{proofsect}[1] {\vskip0.1cm\noindent{\rmfamily\itshape#1.}}{\qed\vspace{0.15cm}}

\theoremstyle{thm}
\newtheorem{theorem}{Theorem}[section]
\newtheorem{lemma}[theorem]{Lemma}
\newtheorem{proposition}[theorem]{Proposition}

\newtheorem*{Main Theorem}{Main Theorem.}

\newtheorem*{special theorem}{Lindeberg-Feller Theorem for Martingales}

\theoremstyle{def}

\theoremstyle{rem}
\newtheorem{remark}[theorem]{{\itshape Remark}}

\numberwithin{equation}{section}

\renewcommand{\section}{\secdef\sct\sect}
\newcommand{\sct}[2][default]{%
\refstepcounter{section}
\addcontentsline{toc}{section}{{\tocsection {}{\thesection}{\!\!\!\!#1\dotfill}}{}}
\vspace{0.7cm}
\centerline{\scshape\thesection.\ #1} \nopagebreak \vspace{0.2cm}}
\newcommand{\sect}[1]{%
\vspace{0.4cm} \centerline{\large\scshape\rmfamily #1}
\vspace{0.2cm}}

\renewcommand{\subsection}{\secdef\subsct\sbsect}
\newcommand{\subsct}[2][default]{\refstepcounter{subsection}
\addcontentsline{toc}{subsection}
{{\tocsection{\!\!}{\hspace{1.2em}\thesubsection}{\!\!\!\!#1\dotfill}}{}}
\nopagebreak\vspace{0.45\baselineskip} {\flushleft\bf
\thesubsection~\bf #1.~}
\\*[3mm]\noindent
\nopagebreak}
\newcommand{\sbsect}[1]{\vspace{0.1cm}\noindent
\textbf{#1.~}\vspace{0.1cm}}

\renewcommand{\subsubsection}{%
\secdef \subsubsect\sbsbsect}
\newcommand{\subsubsect}[2][default]{%
\refstepcounter{subsubsection} 
\addcontentsline{toc}{subsubsection}{{\tocsection{\!\!}
{\hspace{3.05em}\thesubsubsection}{\!\!\!\!#1\dotfill}}{}}
\nopagebreak
\vspace{0.15\baselineskip} \nopagebreak {\flushleft\rmfamily
\itshape\thesubsubsection
\ \rmfamily #1\/.}\ }
\newcommand{\sbsbsect}[1]{\vspace{0.1cm}\noindent
\rmfamily \itshape
\arabic{section}.\arabic{subsection}.\arabic{subsubsection} \
\sffamily #1\/.\ }

\iffinal

\else

\fi

\renewcommand{\caption}[1]{%
\vglue0.5cm
\refstepcounter{figure}
\begin{minipage}{0.9\textwidth}\small {\sc Figure~\thefigure. }#1\end{minipage}}

\newcommand{\textd}{\text{\rm d}\mkern0.5mu}

\newcommand{\texte}{\text{\rm e}}

\newcommand{\1}{\operatorname{\sf 1}}

\newcommand{\mini}{{\text{\rm min}}}

\newcommand{\G}{\mathbb G}

\newcommand{\R}{\mathbb R}

\newcommand{\T}{\mathbb T}

\newcommand{\Z}{\mathbb Z}

\newcommand{\twoeqref}[2]{(\ref{#1}--\ref{#2})}
\newcommand{\cc}{{\text{\rm c}}}

\def\myffrac#1#2 in #3{\raise 2.6pt\hbox{$#3 #1$}\mkern-1.5mu\raise 0.8pt\hbox{$#3/$}\mkern-1.1mu\lower 1.5pt\hbox{$#3 #2$}}
\newcommand{\ffrac}[2]{\mathchoice%
	{\myffrac{#1}{#2} in \scriptstyle}
	{\myffrac{#1}{#2} in \scriptstyle}
	{\myffrac{#1}{#2} in \scriptscriptstyle}
	{\myffrac{#1}{#2} in \scriptscriptstyle}
}

\newcommand{\Bin}{\text{\rm Bin}}
\newcommand{\pt}{p_{\text{\rm\tiny T}}}
\newcommand{\qt}{q_{\text{\rm\tiny T}}}

\title[Bootstrap percolation on trees]
{Metastable behavior for bootstrap\\percolation on regular trees}
\author[M.~Biskup and R.H.~Schonmann]{Marek Biskup$^{1,2}$ \and Roberto H.\ Schonmann$^1$}

\begin{document}
\thanks{\hglue-4.5mm\fontsize{9.6}{9.6}\selectfont\copyright\,2009 by M.~Biskup and R.H.~Schonmann. Reproduction, by any means, of the entire
article for non-commercial purposes is permitted without charge.\vspace{2mm}}
\maketitle

\vspace{-4mm}
\centerline{\textit{$^1$Department of Mathematics, University of California at Los Angeles}}
\centerline{\textit{$^2$School of Economics, University of Southern Bohemia, \v Cesk\'e Bud\v ejovice}}

\vspace{-2mm}
\begin{abstract}
We examine bootstrap percolation on a regular $(b+1)$-ary tree with initial law given by Bernoulli($p$). The sites are updated according to the usual rule: a vacant site becomes occupied if it has at least $\theta$ occupied neighbors, occupied sites remain occupied forever. It is known that, when $b>\theta\geq2$, the limiting density~$q=q(p)$ of occupied sites exhibits a jump at some $\pt=\pt(b,\theta)\in(0,1)$ from $\qt:=q(\pt)<1$ to $q(p)=1$ when~$p>\pt$. We investigate the metastable behavior associated with this transition. Explicitly, we pick $p=\pt+h$ with $h>0$ and show that, as~$h\downarrow0$, the system lingers around the ``critical'' state for time order~$h^{-1/2}$ and then passes to fully occupied state in time~$O(1)$. The law of the entire configuration observed when the occupation density is $q\in(\qt,1)$ converges, as $h\downarrow0$, to a well-defined measure.
\end{abstract}

\bigskip

\section{Introduction}
\noindent
Consider a connected and bounded-degree graph~$\G=(V,E)$. Bootstrap percolation on~$\G$ is a stochastic particle system on the configuration space $\{0,1\}^V$ defined as follows: Given a configuration $\sigma=(\sigma_x)\in\{0,1\}^V$, we interpret $\sigma_x=1$ as ``$x$ is occupied'' and $\sigma_x=0$ as ``$x$ is vacant.''  Pick an integer~$\theta\ge1$ to be referred to as the threshold. Under continuous dynamics, a vacant vertex becomes occupied at rate 1 as soon as it has at least $\theta$ occupied neighbors, otherwise the rate is zero. An occupied vertex remains occupied forever. In discrete-time evolution the updates are performed at all vertices simultaneously 
at regular time intervals (which means that the evolution is entirely deterministic and the model is thus an example of a cellular automaton). The initial configuration is sampled from an i.i.d.\ distribution; i.e., Bernoulli percolation with parameter~$p$. Hence also the name of the process.

Bootstrap percolation arose in physics context in studying the effects of dilution on magnetic properties of crystalline materials (Chalupa, Leath and Reich~\cite{Chalupa-Leath-Reich}). Notwithstanding, the main source of motivation for probabilists was the connection with metastability, and generally non-equilibrium statistical mechanics (Aizenman and Lebowitz~\cite{Aizenman-Lebowitz}). A large body of work exists for the problem on the hypercubic lattice $\G=\Z^d$. In this case, the dynamics converges to fully occupied lattice whenever $p>0$ and $\theta\le d$, and never to fully occupied lattice when $p<1$ and~$\theta>d$ (Schonmann~\cite{Schonmann}). While the ``critical density'' in the most basic case $d=2=\theta$ is zero (van Enter~\cite{vanEnter}), the transition from partial to full occupancy is exhibited in the limit $p\downarrow0$ by the process restricted to boxes of scale $\exp(\ffrac\lambda p)$ as $\lambda$ varies through $\ffrac{\pi^2}{18}$ (see Aizenman and Lebowitz~\cite{Aizenman-Lebowitz}, Holroyd~\cite{Holroyd} 
for results in $d=2$ and those by Cerf and Cirillo~\cite{Cerf-Cirillo}, Cerf and Manzo~\cite{Cerf-Manzo}, Holroyd~\cite{Holroyd2} and Balogh, Bollob\'as and Morris~\cite{Balogh-Bollobas-Morris} dealing with higher dimensions).

In this short note we examine bootstrap percolation for~$\G$ being a homogeneous tree with degree $b+1$. The dynamics, both discrete and continuous-time, can in this case be characterized rather explicitly. (This is not always the case even on a tree; e.g., the stochastic Ising model or kinetically constrained models of Ising type.) It is well known (Chalupa, Leath and Reich~\cite{Chalupa-Leath-Reich}, Balogh, Peres and Pete~\cite{Balogh-Peres-Pete}) that, for $b>\theta\geq2$, there exists a critical density~$\pt$ such that the limiting density of occupied sites converges to $\qt=\qt(p)<1$ when~$p\le\pt$ and to one when~$p>\pt$. There is yet another critical value, $p_\cc$ with $p_\cc<\pt$, above which the terminal configurations will contain infinite connected components (Fontes and Schonmann~\cite{Fontes-Schonmann}).

Our main goal is to point out an interesting metastability effect associated with the first-order transition at~$p=\pt$. Our findings are as follows: For the initial density $p$ slightly above~$\pt$, i.e., $p=\pt+h$, there is a meta\-sta\-ble time scale of length
\begin{equation}
t(h):=\alpha h^{-1/2}+O(1),\qquad h\downarrow0,
\end{equation}
for some~$\alpha\in(0,\infty)$, at which we see a transition from density $\qt$ to full occupancy in a window of times of order unity! This extremely sharp ``cutoff phenomenon'' is perhaps even less intuitive when we add that, as~$h\downarrow0$, the law of the full configuration at the time when density equals some value~$q\in(\qt,1)$ tends to a distinct measure $\mu_q$ (of that density). Thus, a few moments before $t(h)$ we see a distinctly non-fully occupied graph while a few moments after $t(h)$ the graph is essentially full. Presumably, the sharpness of the flip from the metastable to the stable equilibrium is a result of some global event happening at spatial distance --- measured in graph-theoretical sense --- of order $h^{-1/2}$. However, the nature of this event remains mysterious at this point and is definitely worthy of further investigation.

We finish by noting that analogous first order transitions to the one studied here have been found of relevance for jamming (Schwarz, Liu and Chayes~\cite{Lincoln-jamming}). A somewhat related model mimicking a dynamic of single avalanche in directed tree geometry has been thoroughly analyzed by Biskup, Blanchard, Chayes, Gandolfo and Kr\"uger~\cite{Biskup-Chayes-others}.

\section{Recursion equations}
\noindent
We begin by presenting the basic recursion-equation calculations that are permitted by the tree structure of the problem. We will examine two distinct dynamics; first, a discrete time process, when at integer times each unoccupied vertex examines its neighbors, and if there are more than $\theta$ of those occupied, it becomes occupied. The second, continuous time, is such that each vertex is equipped with a clock that rings at exponentially distributed (independent) time intervals, and the update happens accordingly. 

We will consider the process on the $(b+1)$-ary regular tree~$\T_b$, and denote the configuration at time~$t$ by $\sigma(t)=(\sigma_x(t))_{x\in\T_b}$. Here~$\sigma_x(t)\in\{0,1\}$, with~$0$ indicating vacant and~$1$ occupied state, and the initial condition~$\sigma(0)$ is i.i.d.\ with density~$p$. The character of the dynamics implies that~$t\mapsto\sigma_x(t)$ is increasing at each~$x$, and can thus be completely characterized by the quantities
\begin{equation}
T_x:=\inf\{t\ge0\colon\sigma_x(t)=1\bigr\}.
\end{equation}
In particular, $T_x=0$ if~$x$ is occupied initially. Define
\begin{equation}
P_p(t):=P(T_x\le t)
\end{equation}
which, by the symmetries of the dynamics and initial data with respect to tree automorphisms, does not depend on~$x$. Let $\Bin(b,q,\theta)$ denote the probability that a sum of~$b$ i.i.d.\ Bernoulli($q$) is at least~$\theta$, i.e.,
\begin{equation}
\Bin(b,q,\theta):=\sum_{k=\theta}^b\binom bk q^k(1-q)^{b-k}.
\end{equation}
Our first observation is:

\begin{lemma}
\label{lemma-QP}
Let~$Z=1$ for the discrete-time process and $Z=\exp(1)$ for the continuous-time process (independent of the $\sigma$'s). Then
\begin{equation}
\label{P-integral}
P_p(t)=p+(1-p)E\bigl(\Bin(b+1,Q_p(t-Z),\theta)\1_{\{Z\le t\}}\bigr),
\end{equation}
where~$t\mapsto Q_p(t)$ is the unique solution in $[0,1]$ of
\begin{equation}
\label{Q-integral}
Q_p(t)=p+(1-p)E\bigl(\Bin(b,Q_p(t-Z),\theta)\1_{\{Z\le t\}}\bigr).
\end{equation}
In particular, for the continuous time evolution
\begin{equation}
\label{ddtQ}
\frac{\textd}{\textd t}Q_p(t)=W_p\bigl(Q_p(t)\bigr)
\end{equation}
where
\begin{equation}
W_p(q):=p+(1-p)\Bin(b,q,\theta)-q.
\end{equation}
\end{lemma}

\begin{proofsect}{Proof}
Let $\T_b^\star$ denote the rooted tree with forward branching number~$b$ and let~$\varnothing$ denote the root. Consider the bootstrap process restricted to~$\T_b^\star$. We claim that the quantity
\begin{equation}
Q_p(t)\,:=P(T_\varnothing\le t)
\end{equation}
satisfies \eqref{Q-integral}. Indeed, let~$S$ denote the first time the root has at least~$\theta$ occupied neighbors and let~$t-Z$ be the time of the most recent ``clock-ring'' at the root before time~$t$. Then
\begin{equation}
\{0<T_\varnothing\le t\}=\{T_\varnothing>0\}\cap\{S\le T_\varnothing\}\cap\{S\le t-Z\}.
\end{equation}
But on the event $\{S\le T_\varnothing\}$ each of the neighbors at the origin evolves independently  according to same law as the origin itself. Hence, to get the event on the right-hand side, we need that at least~$\theta$ neighbors of the origin are occupied by time~$t-Z$, i.e.,
\begin{equation}
P(0<T_\varnothing\le t)=(1-p)E\bigl(\Bin(b,Q_p(t-Z),\theta)\1_{\{Z\le t\}}\bigr),
\end{equation}
where the prefactor accounts for the probability that, at time zero, the origin was vacant. This immediately yields \eqref{Q-integral}.

Now let us go back to~$\T_b$ and the full bootstrap percolation process on this graph. Let~$x\in\T_b$ be any vertex. If~$S_x$ is the first time~$x$ has at least~$\theta$ occupied neighbors, then $\{0<T_x\le t\}$ can again be written as~$\{T_\varnothing>0\}$ intersected by $\{S_x\le T_\varnothing\}$ and $\{S_x\le t-Z\}$. The difference compared to the ``directed problem'' is that now~$x$ has~$b+1$ neighbors and so
\begin{equation}
P(0<T_x\le t)=(1-p)E\bigl(\Bin(b+1,Q_p(t-Z),\theta)\1_{\{Z\le t\}}\bigr).
\end{equation}
This again immediately gives \eqref{P-integral}.

To get \eqref{ddtQ} we note that, in the case of continuous-time dynamics,~$Z$ is exponentially distributed. A simple change of variables makes \eqref{Q-integral} into
\begin{equation}
Q_p(t)=p+(1-p)\int_0^t\textd z\,\texte^{-(t-z)}\,\Bin\bigl(b,Q_p(z),\theta\bigr)
\end{equation}
This shows that~$t\mapsto Q_p(t)$ is smooth; differentiation then yields \eqref{ddtQ}.
\end{proofsect}

\begin{remark}
For discrete-time evolution, the analogue of \eqref{ddtQ} is
\begin{equation}
Q_p(n+1)-Q_p(n)=W_p\bigl(Q_p(n)\bigr)
\end{equation}
and so, by standard approximation methods, the discrete and continuous problems can be shown to follow similar trajectories. However, to avoid tedious calculations, we will not attempt to cast this statement in quantitative form.
\end{remark}

Since~$Q_p(t)$ readily determines~$P_p(t)$, we will focus the forthcoming discussion on~$Q_p$. The indisputable advantage of continuous time is that \eqref{ddtQ} can be integrated out:
\begin{equation}
t=\int_p^{Q_p(t)}\frac{\textd q}{W_p(q)}.
\end{equation}
The analysis of $Q_p(t)$ thus reduces to the analysis of the denominator $W_p$. Here we notice the following facts (some of which can be traced to calculations in Fontes and Schonmann~\cite{Fontes-Schonmann} and earlier work on this model):

\begin{lemma}
\label{lemma-W}
Suppose $b>\theta\ge2$ and $0<p<1$.
\settowidth{\leftmargini}{(1111)}
\begin{enumerate}
\item[(1)]
There is a $\tilde q\in(0,1)$, independent of~$p$, such that $q\mapsto W_p(q)$ is convex for $q<\tilde q$ and concave for $q>\tilde q$.
\item[(2)]
There is $\tilde p\in(0,1)$ such that $q\mapsto W_p(q)$ is decreasing on $[0,1]$ for $p>\tilde p$ (with a unique local minimum at~$q=1$) while for $p<\tilde p$, there is a $q_\mini(p)\in(0,1)$ such that the local minima of $q\mapsto W_p(q)$ in $[0,1]$ occur at $q\in\{q_\mini(p),1\}$. Both of these are then strict.
\item[(3)]
$p\mapsto W_p(q_\mini(p))$ is increasing, so there is $\pt\in(0,\tilde p)$ such that $W_p(q_\mini(p))<0$ for $p<\pt$ and $W_p(q_\mini(p))>0$ for $p>\pt$.
\item[(4)]
When $W_p(q_\mini(p))\le0$, then $q_\mini(p)>p$.
\end{enumerate}\end{lemma}

\begin{proofsect}{Proof}
To get (1) we notice that, by a coupling of Bernoulli's to independent Uniform([0,1]) random variables, 
\begin{equation}
\frac{\textd}{\textd q}\text{Bin}(b,q,\theta)=b\binom{b-1}{\theta-1}q^{\theta-1}(1-q)^{b-\theta}.
\end{equation}
For~$b>\theta>1$, this starts increasing and then becomes decreasing, implying the said convexity types of $W_p$. Note that $W_p'(0),W_p'(1)<0$ and set $a:=\frac{\textd}{\textd q}\text{Bin}(b,\tilde q,\theta)$. If $W'_p(\tilde q)=a(1-p)-1<0$ then $q\mapsto W_p(q)$ is (strictly) decreasing throughout $[0,1]$, while for $a(1-p)>1$, there is a secondary (strict) local minimum in $(0,\tilde q)$. This yields~(2) with $\tilde p$ defined by $a(1-\tilde p)=1$.

To get~(3), we note that $p\mapsto W_p(q_\mini(p))$ has derivative $1-\text{Bin}(b,q_\mini(p),\theta)$ which is positive for all $0<p<1$. The crossing of zero level must occur as, on the one hand, $W_p(0)=p$ while $W_p'(0)=-1$ and so $q\mapsto W_p(q)$ definitely plunges below zero for small~$p$, while, on the other hand, $W_p$ is strictly positive on $[0,1)$ for~$p$ close to~$\tilde p$. To see~(4), we note that the (strict) convexity on $(0,\tilde q)$ implies 
\begin{equation}
W_p(q)>W_p(0)+W'_p(0)q=p-q
\end{equation}
and so~$W_p(q)>0$ for~$q\le p$.
\end{proofsect}

\begin{remark}
We note that $q\mapsto W_p(q)$ is convex throughout $[0,1]$ when $\theta=1$, while it is concave when $\theta=b$. No local minimum thus appears inside $(0,1)$ for either of these cases and $W_p$ is thus minimized (and vanishes) only at~$q=1$. None of the interesting effects described here occur in either of these cases.
\end{remark}

The findings in Lemma~\ref{lemma-W} imply that that, for~$p<\pt$ --- which by (4) forces~$p<q_\mini(p)$ --- as~$t\to\infty$,~$Q_p(t)$ tends to the smallest positive root of $W_p$ while for~$p>\pt$, the integration eventually passes through the bottleneck at $q_\mini(p)$ under the graph of~$W_p$ and so $Q_p(t)\to1$.  A more detailed look at the bottleneck then yields:

\begin{theorem}
\label{thm-Qscaling}
Suppose $b>\theta\geq2$.
Let~$\pt$ denote the value of~$p$ at which $\min W_p=0$. Set $\qt:=q_\mini(\pt)$ and let $p_h:=\pt+h$. Denote
\begin{equation}
\label{alpha-def}
\alpha:=\frac\pi{\sqrt{\frac12W_{\pt}''(\qt)}}\sqrt{\frac{1-\pt}{1-\qt}}.
\end{equation}
Then for each~$r\in\R$, the limit
\begin{equation}
\label{direct-lim}
\phi(r):=\lim_{h\downarrow0} Q_{p_h}(\alpha h^{-1/2}+r)
\end{equation}
exists and is a solution to the ODE
\begin{equation}
\label{ODE}
\phi'(r)=W_{\pt}\bigl(\phi(r)\bigr),\qquad r\in\R,
\end{equation}
subject to boundary conditions
\begin{equation}
\label{BC}
\lim_{r\to-\infty}\phi(r)=\qt\quad\text{and}\quad\lim_{r\to\infty}\phi(r)=1.
\end{equation}
(This solution is determined uniquely up to a shift of the argument.) 
Equivalently, given~$q\in(\qt,1)$, if~$t_h(q)$ is defined by~$Q_{p_h}(t_h(q))=q$, then
\begin{equation}
\label{th-lim}
t_h(q)-\alpha h^{-1/2}\,\underset{h\downarrow0}\longrightarrow\,\phi^{-1}(q),\qquad \qt<q<1.
\end{equation}
\end{theorem}

This result encompasses the principal findings of this paper. The interpretation is that, for initial density $p$ slightly above~$\pt$, the system evolves into a metastable state close to that for the critical initial density~$\pt$. However, this metastable regime ends rather abruptly at time $\alpha (p-\pt)^{-1/2}+O(1)$ when a sudden transition to full occupancy occurs in a time window of order unity.

\smallskip
\section{Proof of metastability phenomenon}
\noindent
We now commence with the proof of Theorem~\ref{thm-Qscaling}. Key in the calculation is the following observation:

\begin{lemma}
\label{lemma-key}
Let $q<1$ and let~$w$ be a non-negative, $C^3$ function such that
\begin{equation}
w(q)=0=w'(q)\quad\text{and}\quad w''(q)>0.
\end{equation}
Then for all sufficiently small~$\delta>0$,
\begin{equation}
\label{int-asymp}
\lim_{\theta\downarrow0}\,\Bigr[\,\int_{q-\delta}^{q+\delta}\frac{\textd x}{w(x)+\theta(1-x)}\,-\beta\theta^{-1/2}\Bigr]\,\,\,
\text{exists and is finite},
\end{equation}
where $\beta:=\pi/\sqrt{\frac12{w''(q)(1-q)}}$.
\end{lemma}

\begin{proofsect}{Proof}
This follows by standard, albeit tedious, calculations. First, set $x\mapsto q+(1-q)\tilde x$ and $w(x)=(1-q)\tilde w(\tilde x)$ and note that this reduces the problem to~$q=0$. (We also need to note that $w''(q)(1-q)=\tilde w''(0)$ and that $\delta$ gets to be scaled by the factor $(1-q)$.) 

We will henceforth assume that $q=0$ and that $\delta$ is so small that~$w(x)>0$ for all~$x\in[-\delta,\delta]\setminus\{0\}$. First, abbreviate
\begin{equation}
R(x):=\frac1{w(x)+\theta(1-x)}-\frac1{w(x)+\theta}
\end{equation}
and check that
\begin{multline}
\qquad
R(x)+R(-x)
\\=x\theta\frac{[w(-x)-w(x)][w(x)+w(-x)+2\theta+x\theta]}{[w(x)+\theta][w(x)+\theta(1-x)][w(-x)+\theta][w(-x)+\theta(1+x)]}.
\end{multline}
The bounds $w(-x)-w(x)=O(x^3)$ and $x^2/w(x)=O(1)$ show that this expression is bounded uniformly as~$\theta\downarrow0$ and $|x|$~small; the $\theta\downarrow0$ limit is then obviously zero. In light of the $x\leftrightarrow-x$ symmetry of the integration domain and the Bounded Convergence Theorem, we may thus replace the integrand in \eqref{int-asymp} by $[w(x)+\theta]^{-1}$; the claim thus boils down to showing that
\begin{equation}
\lim_{\theta\downarrow0}\,\Bigr[\,\int_{-\delta}^\delta\frac{\textd x}{w(x)+\theta}\,-\beta\theta^{-1/2}\Bigr]\,\,\,
\text{exists and is finite}.
\end{equation}
Here we abbreviate
\begin{equation}
a_2:=\frac12 w''(0)\quad\text{and}\quad a_3:=\frac16w'''(0)
\end{equation}
and invoke the identity
\begin{multline}
\quad
\frac12\frac1{w(x)+\theta}+\frac12\frac1{w(-x)+\theta}-\frac1{a_2x^2+\theta}
\\
=\frac12a_3x^3\frac{w(-x)-w(x)}
{[a_2x^2+\theta][w(x)+\theta][w(-x)+\theta]}\qquad\qquad
\\*[2mm]
-\frac12\frac{w(x)-a_2x^2-a_3x^3}{[w(x)+\theta][a_2x^2+\theta]}-\frac12\frac{w(-x)-a_2x^2+a_3x^3}{[w(-x)+\theta][a_2x^2+\theta]}
\quad
\end{multline}
Again, all three expressions on the right-hand side are bounded uniformly in $x\in[-\delta,\delta]$  and have a well-defined limit as~$\theta\downarrow0$. The Bounded Convergence Theorem again ensures that the integral of $[w(x)+\theta]^{-1}$ equals that of $[a_2x^2+\theta]^{-1}$ plus a quantity that has a finite limit as~$\theta\downarrow0$. The claim now follows by noting that
\begin{equation}
\int_{-\delta}^\delta\frac{\textd x}{a_2x^2+\theta}=\frac\pi{\sqrt{a_2\theta}}-\frac2\delta+o(1),
\qquad\theta\downarrow0,
\end{equation}
which is checked by a direct calculation.
\end{proofsect}

\begin{proofsect}{Proof of Theorem~\ref{thm-Qscaling}}
It suffices to prove \eqref{th-lim}; indeed, the claim \eqref{direct-lim} then follows by the strict monotonicity of~$t\mapsto Q_p(t)$ and of the solution to the ODE \eqref{ODE} subject to \eqref{BC}. First note that
\begin{equation}
t_h(q)=\int_{p_h}^q\frac{\textd s}{W_{p_h}(s)}.
\end{equation}
Thus, the point-wise limit $W_{p_h}(s)\to W_{\pt}(s)$ and the fact that $W_{\pt}(s)>0$ for~$s>\qt$ imply that $t_h(q)-t_h(q')$, for $q,q'>\qt$, converges to $\phi^{-1}(q)-\phi^{-1}(q')$, for some solution~$r\mapsto\phi(r)$ to \twoeqref{ODE}{BC}. To show \eqref{th-lim} it now suffices to prove that (for at least one~$q\in\R$) 
\begin{equation}
\label{th-2ndlim}
\lim_{h\downarrow0}\bigl[\,t_h(q)-\alpha h^{-1/2}\bigr]\quad\text{exists and is finite}.
\end{equation}
A calculation yields
\begin{equation}
W_{p_h}(s)=\frac{1-p_h}{1-\pt}\left(W_{\pt}(s)+\theta_h(1-s)\right)
\end{equation}
with $\theta_h:=h/(1-p_h)$. Noting that $W_{\pt}$ is uniformly positive on the complement of the interval $(\qt-\delta,\qt+\delta)$ and observing the asymptotics
\begin{equation}
\frac{1-p_h}{1-\pt}-1=o(h^{1/2})
\quad\text{and}\quad
\frac{\theta_h}h-(1-\pt)=o(h^{1/2}),
\end{equation}
to get \eqref{th-2ndlim} it thus suffices to show that 
\begin{equation}
\label{lim-exists}
\lim_{\theta\downarrow0}\Bigl[\,\int_{\qt-\delta}^{\qt+\delta}\frac{\textd s}{W_{\pt}(s)+\theta(1-s)}-\beta\theta^{-1/2}\Bigr]\quad\text{exists and is finite}
\end{equation}
for~$\beta:=\pi/\sqrt{\frac12W_{\pt}''(\qt)(1-\qt)}$. But this follows from Lemma~\ref{lemma-key}.
\end{proofsect}

\begin{remark}
The asymptotic analysis in the above proofs allows us to control the speed of convergence at the times way before (or after) $\alpha h^{-1/2}+O(1)$. Thus one can prove, for instance, that for all $0<\lambda<\alpha$,
\begin{equation}
h^{-1/2}\bigr[
Q_{p_h}(\lambda h^{-1/2})-\qt
\bigr]
\,\underset{h\downarrow0}\longrightarrow c_1\tan\bigl(c_2\lambda-c_3\bigr)
\end{equation}
for some constants~$c_1,c_2>0$ and $c_3\in\R$.
\end{remark}

\section{Process convergence}
\label{sec3}
\noindent
The above calculations give us information about the evolution of the law at a single vertex. However, as we will show in this section, the tree structure and some additional arguments allow us to control the law of the entire infinite-volume configuration at times close to and after passing through the metastable regime. For~$t\ge0$, let $\mu_{p,t}$ to denote the law of $(\sigma_x(t))$ started from Bernoulli$(p)$.

\begin{theorem}
\label{thm3}
There is a stochastically increasing family of tree-autormorphism invariant measures~$(\nu_q)_{0\le q\le1}$ on $\{0,1\}^{\T_b}$ with density $\nu_q(\sigma_x=1)=\text{\rm Bin}(b+1,q,\theta)$ such that:
\settowidth{\leftmargini}{(1111)}
\begin{enumerate}
\item[(1)]
If~$p\le \pt$ and $q:=\lim_{t\to\infty}Q_p(t)$, then
\begin{equation}
\mu_{p,t}\overset{\text{\rm w}}\longrightarrow\nu_q,\qquad t\to\infty.
\end{equation}
\item[(2)]
If~$p>\pt$ then
\begin{equation}
\mu_{p,t}\overset{\text{\rm w}}\longrightarrow\nu_1:=\delta_{\underline1},\qquad t\to\infty.
\end{equation}
\end{enumerate}
In addition, we also have:
\begin{enumerate}
\item[(3)]
Set $p_h=\pt+h$ and, given $q\in(\qt,1)$, let $t_h(q)$ be defined by $Q_{p_h}(t_h(q))=q$. Then
\begin{equation}
\mu_{p_h,t_h(q)}\overset{\text{\rm w}}\longrightarrow\nu_q,\qquad h\downarrow0.
\end{equation}
\end{enumerate}
\end{theorem}

\begin{remark}
We note that the statement (3) is valid only for the continuous-time dynamics as, for discrete time, even the occupation density changes in discrete quanta.
\end{remark}

The convergence in (1-2) follows trivially from the fact that $t\mapsto\mu_{p,t}$ is stochastically increasing; the stochastic ordering of $q\mapsto\nu_q$ is a consequence of the existence of monotone coupling of the processes with different initial densities and --- in part (3) --- also the monotonicity of the dynamics and of the function $q\mapsto t_h(q)$. The main issue is thus the proof of (3). The key observation is the fact that the effect of the rest of the system on a finite set can be represented through a time-dependent boundary condition. 

\smallskip
Let~$A$ be a finite connected subset of~$\T_b$ and let~$\partial A$ denote the vertices in~$\T_b\setminus A$ that have a neighbor in~$A$. For each~$x\in\partial A$ let~$\T(x)$ denote the (unique) subtree of~$\T_b\setminus A$ rooted in~$x$. For any initial condition $\sigma_x(0)=\sigma_x$, let $\{\sigma_x(t)\colon x\in\T_b\}$ be the state of bootstrap percolation on~$\T_b$ at time~$t$. In addition, for each~$x\in\partial A$, let $\{\eta_{A,z}(t)\colon z\in\T(x)\}$ be the state at time~$t$ of bootstrap percolation on the graph~$\T(x)$ which is started from the same initial condition as~$\sigma$. Given a trajectory $\{\eta_{A,z}(t)\colon z\in\partial A,\,t\ge0\}$, we now define $\{\eta_{A,x}(t)\colon x\in A\}$ to be the state of bootstrap percolation on~$A$ at time~$t$ that initiates from $\eta_{A,x}=\sigma_x$ and evolves against the time-dependent ``boundary condition'' $\eta_{A,z}(t)$ at all~$z\in\partial A$. 

\begin{lemma}
\label{lemma3.3}
The laws of $\{\sigma_x(t)\colon x\in A,\,t\ge0\}$ and $\{\eta_{A,x}(t)\colon x\in A,\,t\ge0\}$ are identical.
\end{lemma}

\begin{proofsect}{Proof}
For each~$x\in\partial A$, let $n(x)\in A$ be the unique neighbor of~$x$ in~$A$. We can now couple the evolutions of $\{\sigma_x(t)\colon x\in\T_b,\,t\ge0\}$ and $\{\eta_{A,x}(t)\colon x\in\T_b,\,t\ge0\}$ by using the same rings of Poisson clocks. Given a sample of all these rings on~$\T_b$, for all~$x\in\partial A$ we have
\begin{equation}
\sigma_z(t)=\eta_{A,z}(t),\qquad z\in\T(x),\,0\le t<T_{n(x)}.
\end{equation}
However, at~$t=T_{n(x)}$ the vertex~$n(x)$ becomes occupied and, since~$\T(x)$ has no other neighbor in~$A$, the evolution on~$\T(x)$ becomes irrelevant for that in~$A$. It follows that $\sigma_z(t)=\eta_{A,z}(t)$ for all~$z\in A$ and all~$t\ge0$ which implies the claim.
\end{proofsect}

As we will see, Theorem~\ref{thm3} will be a consequence of this, slightly stronger, claim:

\begin{proposition}
\label{prop4}
Let $\alpha$ be as in \eqref{alpha-def} and let $h\mapsto r_h$ be such that both $r_h$ and $\alpha h^{-1/2}-r_h$ tend to infinity as~$h\downarrow0$. Define $S_x\in\{0,1\}\times\R$ via
\begin{equation}
\label{3.5}
S_x:=\bigl(\1_{\{T_x\le r_h\}}, T_x-\alpha h^{-1/2}\1_{\{T_x> r_h\}}\bigr).
\end{equation}
As~$h\downarrow0$, the law of $(S_x)$ induced by $P_p$ with~$p=\pt+h$ on the product space $(\{0,1\}\times\R)^{\T_b}$ tends to a probability measure which is independent of the choice of~$r_h$. 
\end{proposition}

\begin{proofsect}{Proof}
We begin by noting that the probability of the event
\begin{equation}
\label{event}
\{T_x\le r_h\}\cup\{\alpha h^{-1/2}-r_h\le T_x\le\alpha h^{-1/2}+r_h\}
\end{equation}
tends to one as~$h\downarrow0$, for any~$r_h$ as above. Indeed, the events are disjoint once $h$ is sufficiently small, the probability of the former is at least $\qt$ in the $h\downarrow0$ limit by the fact that $Q_{\pt+h}(t)\ge Q_{\pt}(t)\to\qt$ while the probability of the latter tends to $1-\qt$ by \eqref{th-lim}. It follows that the cutoff to $T_x>r_h$ may be replaced by
\begin{equation}
\alpha h^{-1/2}-r_h\le T_x\le\alpha h^{-1/2}+r_h
\end{equation}
in the second half of \eqref{3.5}. 

To describe the convergence of the full distribution, pick a finite set~$A\subset\T_b$ and, given the initial data $\sigma_A:=\{\sigma_x\colon x\in A\}$ and boundary data $S_{\partial A}:=\{S_x\colon x\in\partial A\}$, let $\mu_A(-|S_{\partial A},\sigma_A)$ be the bootstrap percolation process on~$A$ started from the initial configuration~$\sigma_A$ and evolving against the time-dependent boundary data~$S_{\partial A}$. Lemma~\ref{lemma3.3} tells us that the law of the variables $\{S_x\colon x\in A\}$ is that of $\mu_A(-|S_{\partial A},\sigma_A)$ integrated over Bernoulli($p$) variables~$\sigma_A$ and i.i.d.\ boundary data where each~$S_x$, $x\in\partial A$, is sampled independently from the law of~$S_\varnothing$ at the root of a rooted tree. 

Recall that, by \eqref{th-lim} the law of $S_{\partial A}$ has a non-degenerate weak limit as~$h\downarrow0$, which is the same regardless of how~$r_h\to\infty$ and $\alpha h^{-1/2}-r_h\to\infty$. By definition, also the law of $\sigma_A$ trivially converges as $h\downarrow0$. As there is only a finite number of possible $\sigma_A$'s, to prove convergence it suffices to show that, for each $\sigma_A$, the family of measure-valued functions
\begin{equation}
\label{MVF}
S_{\partial A}\mapsto\mu_A(-|S_{\partial A},\sigma_A)
\end{equation}
is equicontinuous. For this we note that a change in $T_x$ for~$x\in\partial A$ by~$\epsilon>0$ will only be felt inside~$A$ if the Poisson clock at the (unique) neighbor~$y\in A$ of~$x$ clicks within the interval $[T_x,T_x+\epsilon]$. This has probability of order~$\epsilon$. The functions \eqref{MVF} are thus Lipschitz continuous in the boundary values and convergence follows.
\end{proofsect}

\begin{proofsect}{Proof of Theorem~\ref{thm3}(3)}
Assuming $r_h\le t_h$, we have
\begin{equation}
\sigma_x(t_h)=\1_{\{T_x\le r_h\}}+\1_{\{T_x>r_h\}}\1_{\{T_x\le t_h\}}.
\end{equation}
The convergence of $\mu_{p_h,t_h}$ then follows from Proposition~\ref{prop4} and \eqref{th-lim}.
\end{proofsect}


\section*{Acknowledgments}
\noindent 
The research of M.B.~was partially supported by the grants NSF DMS-0505356 and NSF DMS-0806198. 
The research of R.H.S.~was partially supported by the grant NSF DMS-0300672.

\end{document}